\documentclass[preprint]{elsarticle}

\usepackage{amsmath,amssymb,amsthm}

\journal{Fuzzy Sets and Systems}

\begin{document}

\begin{frontmatter}




\title
{$L$-fuzzy strongest postcondition predicate transformers as
$L$-idempotent linear or affine operators between semimodules of
monotonic predicates}

\tnotetext[supp]{This research was supported by the~Slovenian Research
Agency grants P1-0292-0101, J1-4144-0101.}


\author[nyk]{Oleh Nykyforchyn\corref{cor1}}
\ead{oleh.nyk@gmail.com}
\author[rep]{Du\v san Repov\v s}
\ead{dusan.repovs@guest.arnes.si}
\cortext[cor1]{}

\address[nyk]{Vasyl' Stefanyk Precarpathian National University,
Shevchenka 57, Ivano-Frankivsk, 76025, Ukraine}
\address[rep]{Faculty of Education and Faculty of Mathematics and Physics,
University of Ljubljana, Jadranska 19, Ljubljana, 1000, Slovenia}

\begin{abstract}
For a~completely distributive quantale $L$, $L$-fuzzy strongest
postcondition predicate transformers are introduced, and it is
shown that, under reasonable assumptions, they are linear or affine
continuous mappings between continuous $L$-idempotent semimodules
of $L$-fuzzy monotonic predicates.
\end{abstract}


\end{frontmatter}

\def\baselinestretch{1.1}

\let\phi\varphi
\let\kappa\varkappa
\let\le\leqslant
\let\ge\geqslant
\let\de\delta
\let\eps\varepsilon

\let\emptyset\varnothing
\let\hra\hookrightarrow
\let\ES\varnothing
\let\subsetne\subsetneq
\newcommand\subcl{\mathrel{\underset{\mathrm{cl}}{\subset}}}
\newcommand\subop{\mathrel{\underset{\mathrm{op}}{\subset}}}

\newcommand\cns{\mathop{\mathrm{const}}}
\newcommand\Int{\mathop{\mathrm{Int}}}
\newcommand\Cl{\mathop{\mathrm{Cl}}}
\newcommand{\cv}{\mathop{\mathrm{conv}}}
\newcommand{\cc}{\mathop{\mathrm{cc}}}
\newcommand{\supp}{\mathop{\mathrm{supp}}}
\newcommand{\supf}{\mathop{\mathrm{supp}_F}}
\newcommand{\bohr}{\mathop{\mathit{bohr}}}
\newcommand{\olim}{\mathop{\overline{\mathrm{lim}}}\limits}
\newcommand{\psup}{\mathop{\mathrm{sup}_p}\limits}
\newcommand{\pinf}{\mathop{\mathrm{inf}_p}\limits}

\newcommand{\pr}{\mathop{\mathrm{pr}}\nolimits}
\newcommand{\rel}{\mathop{\mathrm{rel}}\nolimits}
\newcommand{\diam}{\mathop{\mathrm{diam}}}
\newcommand{\cov}{\mathop{\mathit{cov}}}
\newcommand{\uni}[1]{\mathbf{1}_{{#1}}}
\newcommand{\id}[1]{\mathbf{1}_{#1}}

\newcommand{\sub}{\mathop{\mathrm{sub}}}
\newcommand{\msub}{\mathop{\mathrm{msub}}}
\newcommand{\epi}{\mathop{\mathrm{epi}}}

\newcommand{\expl}{\mathop{\mathrm{exp}_{\,l}}}
\newcommand{\expu}{\mathop{\mathrm{exp}_{\,u}}}
\newcommand{\expus}{\mathop{\mathrm{exp}_{\,us}}}
\newcommand{\explus}{\mathop{\mathrm{exp}^{\le}_{\,us}}}

\newcommand{\sast}{\{\ast\}}
\newcommand{\sastu}{\{\ast,e\}}
\newcommand{\sint}{\sideset{}{^\lor}\int\limits}
\newcommand{\jint}{\sideset{}{^\lor}\int\limits}

\newcommand{\al}{A_{\alpha}}
\newcommand{\alup}{A^{\alpha}}
\newcommand{\ul}{U^{\alpha}}

\newcommand{\bOd}{\mathop{\bar{O}_\delta}}
\newcommand{\bOdd}{\mathop{\bar{O}_{2\delta}}}
\newcommand{\bBd}{\mathop{\bar{B}_\delta}}
\newcommand{\bBdd}{\mathop{\bar{B}_{2\delta}}}
\newcommand{\Mup}{\mathbb{M}_{\scriptscriptstyle\cap}}
\newcommand{\Mdn}{\mathbb{M}_{\scriptscriptstyle\cup}}
\newcommand{\mup}{M_{\scriptscriptstyle\cap}}
\newcommand{\mdn}{M_{\scriptscriptstyle\cup}}
\newcommand{\etaup}{\eta_{\scriptscriptstyle\cap}}
\newcommand{\muup}{\mu_{\scriptscriptstyle\cap}}
\newcommand{\etadn}{\eta_{\scriptscriptstyle\cup}}
\newcommand{\mudn}{\mu_{\scriptscriptstyle\cup}}
\newcommand{\xidn}{\xi_{\scriptscriptstyle\cup}}
\newcommand{\xiup}{\xi_{\scriptscriptstyle\cap}}
\newcommand{\mxup}{M_{\scriptscriptstyle\cap}X}
\newcommand{\mxdn}{M_{\scriptscriptstyle\cup}X}
\newcommand{\ml}{M_{L}}
\newcommand{\Ml}{\mathbb{M}_{L}}
\newcommand{\mldn}{M_{{\scriptscriptstyle\cup} L}}
\newcommand{\Mldn}{\mathbb{M}_{{\scriptscriptstyle\cup} L}}
\newcommand{\mlup}{M_{{\scriptscriptstyle\cap} L}}
\newcommand{\Mlup}{\mathbb{M}_{{\scriptscriptstyle\cap} L}}
\newcommand{\etalup}{\eta_{{\scriptscriptstyle\cap}L}}
\newcommand{\mulup}{\mu_{{\scriptscriptstyle\cap}L}}
\newcommand{\etaldn}{\eta_{{\scriptscriptstyle\cup}L}}
\newcommand{\muldn}{\mu_{{\scriptscriptstyle\cup}L}}
\newcommand{\Delom}{\Delta_{\scriptscriptstyle\otimes}}
\newcommand{\Delop}{\Delta_{\scriptscriptstyle\oplus}}
\newcommand{\Delol}{\Delta_{\scriptscriptstyle L}}
\newcommand{\Delot}{\Delta_{\scriptscriptstyle \tilde L}}
\newcommand{\op}{\oplus}
\newcommand{\om}{\otimes}
\newcommand{\od}{\odot}
\newcommand{\os}{\mathbin{\bar\ast}}
\newcommand{\bp}{\mathop{\bar\oplus}}
\newcommand{\bm}{\mathop{\bar\otimes}}
\newcommand{\bd}{\mathop{\bar\odot}}
\newcommand{\bt}{\mathbin{\underset{{\displaystyle\bar{}}}\odot}}
\newcommand{\bs}{\mathop{\bar\circledast}}
\newcommand\uM{{\rlap{\underline{\phantom{N}}}M}}
\newcommand\oM{{M\llap{$\overline{\phantom{S}}$\,}}}
\newcommand\uO{{\rlap{\underline{\phantom{I\,\,}}}O}}
\newcommand\oO{{O\llap{$\overline{\phantom{I}}$\,}}}
\newcommand\uS{{\underline{\mathcal{S}}}}
\newcommand\oS{{\overline{\mathcal{S}}}}
\newcommand\uP{{\underline{P}}}
\newcommand\oP{{\overline{P}}}
\newcommand\usp{{\rlap{\underline{\phantom{ni}}}{sp}}}
\newcommand\RML{\mathbf{M}_{L}}
\newcommand\tRML{\tilde{\mathbf{M}}_{L}}
\newcommand\RMtL{\mathbf{M}_{\tilde L}}
\newcommand\tRMtL{\tilde{\mathbf{M}}_{\tilde L}}
\newcommand{\Fle}{F^{\scriptscriptstyle \le}_L}
\newcommand{\etale}{\eta^{\scriptscriptstyle \le}_L}
\newcommand{\mule}{\mu^{\scriptscriptstyle \le}_L}
\newcommand{\FF}{\bar F}
\newcommand{\BTT}{\tilde{\mathbb{T}}}
\newcommand{\TT}{\tilde T}
\newcommand{\HT}{\hat T}
\newcommand\hhat[1]{\Hat{\hat #1}}
\newcommand{\tG}{\tilde G}
\newcommand{\tbG}{\tilde{\mathbb{G}}}
\newcommand{\mmu}{\tilde\mu}
\newcommand{\eeta}{\tilde\eta}
\newcommand{\tlp}{\mathbin{\tilde\oplus}}
\newcommand{\tom}{\mathbin{\tilde\otimes}}
\newcommand{\tod}{\mathbin{\tilde\odot}}
\newcommand\tle{\mathrel{\tilde{\le}}}
\newcommand\tge{\mathrel{\tilde{\ge}}}
\newcommand\tland{\mathbin{\tilde{\land}}}
\newcommand\tlor{\mathbin{\tilde{\lor}}}
\newcommand\tsup{\mathop{\tilde{\sup}}}
\newcommand\tinf{\mathop{\tilde{\smash[t]{\inf}}}}
\newcommand\ups{{\uparrow}}
\newcommand\dns{{\downarrow}}
\newcommand\tups{{\tilde\uparrow}}
\newcommand\tdns{{\tilde\downarrow}}
\newcommand{\Exp}{\mathop{\mathrm{Exp}}}
\newcommand\exps{\exp_{\triangle}}
\newcommand\etas{\eta_{\triangle}}
\newcommand\mus{\mu_{\triangle}}
\newcommand\expdd{\exp_{\dns}}
\newcommand\etad{\eta_{\dns}}
\newcommand\mud{\mu_{\dns}}
\newcommand\expuu{\exp_{\ups}}
\newcommand\etau{\eta_{\ups}}
\newcommand\muu{\mu_{\ups}}
\newcommand\expre{\exp_{\scriptscriptstyle\supset}}
\newcommand\etare{\eta_{\scriptscriptstyle\supset}}
\newcommand\expL{\exp^L_\triangle}
\newcommand\etaL{\eta^L_\triangle}
\newcommand\mua{\mu^{*}}
\newcommand\expul{\exp^L_{\ups}}
\newcommand\etaul{\eta^L_{\ups}}
\newcommand\muua{\mu^{*}_{\ups}}
\newcommand\expdl{\exp^L_{\dns}}
\newcommand\etadl{\eta^L_{\dns}}
\newcommand\muda{\mu^{*}_{\dns}}
\newcommand\expLb{\bar{\exp}^L_\triangle}
\newcommand\etaLb{\bar{\eta}^L_\triangle}
\newcommand\muab{\bar{\mu}^{*}_{\triangle}}
\newcommand\expulb{\bar{\exp}^L_{\ups}}
\newcommand\etaulb{\bar{\eta}^L_{\ups}}
\newcommand\muuab{\bar{\mu}^{*}_{\ups}}
\newcommand\expdlb{\bar{\exp}^L_{\dns}}
\newcommand\etadlb{\bar{\eta}^L_{\dns}}
\newcommand\mudab{\bar{\mu}^{*}_{\dns}}
\newcommand\Uland{U_{\land}}
\newcommand\Us{U_{\triangle}}
\newcommand\Hs{H_{\triangle}}
\newcommand\Ud{U_{\dns}}
\newcommand\Hd{H_{\dns}}
\newcommand\Uu{U_{\ups}}
\newcommand\Hu{H_{\ups}}
\newcommand\Ua{U^{*}}
\newcommand\Uua{U^{*}_{\ups}}
\newcommand\Uda{U^{*}_{\dns}}
\newcommand\Usb{\bar{U}_{\triangle}}
\newcommand\Hsb{\bar{H}_{\triangle}}
\newcommand\Udb{\bar{U}_{\dns}}
\newcommand\Hdb{\bar{H}_{\dns}}
\newcommand\Uub{\bar{U}_{\ups}}
\newcommand\Hub{\bar{H}_{\ups}}
\newcommand\Uab{\bar{U}^{*}}
\newcommand\Uuab{\bar{U}^{*}_{\ups}}
\newcommand\Udab{\bar{U}^{*}_{\dns}}

\newcommand{\lip}{\mathrm{Lip}}

\newcommand\Ob{\mathop{\mathrm{Ob}}}
\newcommand\Ar{\mathop{\mathrm{Ar}}}
\newcommand\dom{\mathop{\mathrm{dom}}}
\newcommand\rng{\mathop{\mathrm{rng}}}
\let\ccirc\circledcirc
\newcommand\adot{\mathbin{\underset{\raisebox{.25ex}{\rm *}}{\odot}}}
\newcommand\acirc{\mathbin{\underset{\raisebox{.2ex}{\rm *}}{\circledcirc}}}
\newcommand\bcirc{\mathbin{\bar\circledcirc}}
\newcommand\bacirc{\mathbin{\underset{\raisebox{.2ex}{\rm *}}{\bar\circledcirc}}}
\newcommand\toK{\to_{\CCK}}
\newcommand\ccirK{\mathop{\circledcirc_{\CCK}}}
\newcommand\cirK{\mathop{\circ_{\CCK}}}
\newcommand{\Lip}{\mathop{\rm Lip}}
\newcommand{\Set}{\mathcal{S}\mathrm{et}}
\newcommand{\Top}{{\rlap{\mathsurround=0pt $\mathcal{T}$}\hphantom{\mathrm{n}}\mathrm{op}}}
\newcommand{\Metr}{{\mathcal{M}\mathrm{etr}}}
\newcommand{\Comp}{{\rlap{\mathsurround=0pt $\mathcal{C}$}\hphantom{\mathrm{n}}\mathrm{omp}}}
\newcommand{\Tych}{{\rlap{\mathsurround=0pt $\mathcal{T}$}\hphantom{\mathrm{n}}\mathrm{ych}}}
\newcommand{\Csgr}{\mathcal{CS}\mathrm{gr}}
\newcommand{\Cmon}{\mathcal{C}\mathrm{mon}}
\newcommand{\Cmsgr}{\mathcal{CMS}\mathrm{gr}}
\newcommand{\Casgr}{\mathcal{CA}\mathrm{b}\mathcal{S}\mathrm{gr}}
\newcommand{\Cam}{\mathcal{CA}\mathrm{b}\mathcal{M}\mathrm{on}}
\newcommand{\Conv}{\mathcal{C}\mathrm{onv}}
\newcommand{\Convmm}{\mathcal{C}\mathrm{onv}_{\scriptstyle\max,\min}}
\newcommand{\BConvmm}{\mathcal{B}\mathrm{i}\mathcal{C}\mathrm{onv}_{\scriptstyle\max,\min}}
\newcommand{\Convl}{\mathcal{C}\mathrm{onv}_{L}}
\newcommand{\BConvl}{\mathcal{B}\mathrm{i}\mathcal{C}\mathrm{onv}_{L}}
\newcommand{\Sconv}{\mathcal{SC}\mathrm{onv}}
\newcommand{\Ssconv}{\mathcal{SSC}\mathrm{onv}}
\newcommand{\Wsconv}{\mathcal{WSC}\mathrm{onv}}
\newcommand{\Cla}{\mathcal{CL}}
\newcommand{\Cld}{\mathcal{CL}_{\dns}}
\newcommand{\Clu}{\mathcal{CL}_{\ups}}
\newcommand{\Csem}{\mathcal{CS}\mathrm{em}}
\newcommand{\Cscsem}{\mathcal{CSCS}\mathrm{em}}
\newcommand\Laws{\mathcal{L}\mathrm{aws}}
\newcommand\Llaws{\mathcal{LL}\mathrm{aws}}
\newcommand\Lawsd{\mathcal{L}\mathrm{aws}_{\dns}}
\newcommand\Llawsd{\mathcal{LL}\mathrm{aws}_{\dns}}
\newcommand\Lawsu{\mathcal{L}\mathrm{aws}_{\ups}}
\newcommand\Llawsu{\mathcal{LL}\mathrm{aws}_{\ups}}
\newcommand\Lwsm{(L,\op,*)\text{-}\mathcal{L}\mathrm{w}\mathcal{SM}\mathrm{od}}
\newcommand\Lwsmu{(L,\op,*)\text{-}\mathcal{L}\mathrm{w}\mathcal{SM}\mathrm{od}_{\ups}}
\newcommand\Lwsmd{(L,\op,*)\text{-}\mathcal{L}\mathrm{w}\mathcal{SM}\mathrm{od}_{\dns}}
\newcommand\Lwsa{(L,\op,*)\text{-}\mathcal{L}\mathrm{w}\mathcal{SA}\mathrm{ff}}
\newcommand\Lwsau{(L,\op,*)\text{-}\mathcal{L}\mathrm{w}\mathcal{SA}\mathrm{ff}_{\ups}}
\newcommand\Lwsad{(L,\op,*)\text{-}\mathcal{L}\mathrm{w}\mathcal{SA}\mathrm{ff}_{\dns}}
\newcommand\Camb{\mathcal{CA}\mathrm{mb}}
\newcommand\Csamb{\mathcal{CSA}\mathrm{mb}}
\newcommand\Cpamb{\mathcal{CPA}\mathrm{mb}}
\newcommand\Cpsamb{\mathcal{CPSA}\mathrm{mb}}
\newcommand\Coamb{\mathcal{COA}\mathrm{mb}}
\let\sms\smallsmile
\newcommand\ssms{{%
\rlap{$\scriptstyle\smallsmile$}%
\raise .5ex \hbox{$\scriptstyle\smallsmile$}%
}}

\newcommand{\BBN}{\mathbb{N}}
\newcommand{\BBZ}{\mathbb{Z}}
\newcommand{\BBQ}{\mathbb{Q}}
\newcommand{\BBR}{\mathbb{R}}
\newcommand{\BBC}{\mathbb{C}}
\newcommand{\BBF}{\mathbb{F}}
\newcommand{\BBH}{\mathbb{H}}
\newcommand{\BBI}{\mathbb{I}}
\newcommand{\BBG}{\mathbb{G}}
\newcommand{\BBP}{\mathbb{P}}
\newcommand{\BBK}{\mathbb{K}}
\newcommand{\BBL}{\mathbb{L}}
\newcommand{\BBM}{\mathbb{M}}
\newcommand{\BBT}{\mathbb{T}}
\newcommand{\CCA}{\mathcal{A}}
\newcommand{\CCB}{\mathcal{B}}
\newcommand{\CCC}{\mathcal{C}}
\newcommand{\CCD}{\mathcal{D}}
\newcommand{\CCE}{\mathcal{E}}
\newcommand{\CCF}{\mathcal{F}}
\newcommand{\CCG}{\mathcal{G}}
\newcommand{\CCH}{\mathcal{H}}
\newcommand{\CCI}{\mathcal{I}}
\newcommand{\CCJ}{\mathcal{J}}
\newcommand{\CCK}{\mathcal{K}}
\newcommand{\CCL}{\mathcal{L}}
\newcommand{\CCM}{\mathcal{M}}
\newcommand{\CCN}{\mathcal{N}}
\newcommand{\CCR}{\mathcal{R}}
\newcommand{\CCP}{\mathcal{P}}
\newcommand{\CCS}{\mathcal{S}}
\newcommand{\CCU}{\mathcal{U}}
\newcommand{\CCT}{\mathcal{T}}
\newcommand{\CCV}{\mathcal{V}}
\newcommand{\CCW}{\mathcal{W}}
\newcommand{\RRA}{\mathrm{A}}
\newcommand{\RRB}{\mathrm{B}}
\newcommand{\RRC}{\mathrm{C}}
\newcommand{\RRD}{\mathrm{D}}
\newcommand{\RRE}{\mathrm{E}}
\newcommand{\RRF}{\mathrm{F}}
\newcommand{\RRG}{\mathrm{G}}
\newcommand{\RRI}{\mathrm{I}}
\newcommand{\RRJ}{\mathrm{J}}
\newcommand{\RRK}{\mathrm{K}}
\newcommand{\RRL}{\mathrm{L}}
\newcommand{\RRN}{\mathrm{N}}
\newcommand{\RRP}{\mathrm{P}}
\newcommand\SSA{\mathsf{A}}
\newcommand\SSB{\mathsf{B}}
\newcommand\SSC{\mathsf{C}}
\newcommand\SSD{\mathsf{D}}
\newcommand\SSF{\mathsf{F}}


\newtheorem{theorem}{Theorem}
\newtheorem{proposition}{Proposition}[section]
\newtheorem{lemma}[proposition]{Lemma}
\newtheorem{corollary}[proposition]{Corollary}
\newtheorem{definition}[proposition]{Definition}

\theoremstyle{remark}
\newtheorem{example}[proposition]{Example}
\newtheorem{remark}[proposition]{Remark}
\newtheorem*{remark*}{Remark}

\section*{Introduction}

Predicate transformers, which were introduced in the~pioneering
work of Dijkstra~\cite{Dijk:Guarded:75}, are powerful tools for
analyzing the~total or partial correctness of computer programs.
The~main idea is that a~final state after execution of a~program
depends on its initial state; hence there is an~interdependency
between validity of statements (predicates) about the~initial and
the~final states. One can ask, e.g., what are minimal requirements
on an~initial state that ensure that the~final state satisfies
a~certain condition. Then these requirements form the~weakest
precondition for the~given condition. On the~other hand, the~most
precise knowledge about an~output of a~program for an~input, that
satisfies some predicate, is the~strongest postcondition for this
predicate. Such ``forward'' and ``backward'' dependencies are
called predicate transformers.

Things become more complicated because of randomness or/and
non-de\-ter\-mi\-nism, which can arise from unpredictable
influence, ``angelic'' or ``demonic'' (with the~obvious
connotations). For simplicity, assume first that only randomness is
present, and a~set $S$ of possible states is finite. We mostly
follow~\cite{MMS:ProbPredTrans:96}, but notation will partially
vary. A~\emph{subprobabilistic distribution} $D:S\to [0,1]$
guarantees that the~probability of each state $s\in S$ is at
least~$D(s)$. Obviously it is required that $\sum_{s\in S}D(s)\le
1$, and $1-\sum_{s\in S}D(s)$ ``goes to'' unspecified state of
the~system. We say that a~subprobabilistic distribution $D$ is
\emph{refined} by another subprobabilistic distribution $D'$ on $S$
(written $D\sqsubseteq D'$) if $D(s)\le D'(s)$ for all $s\in S$;
this means that $D'$ offers more precise knowledge than~$D$. This
partial order makes the~set $\bar S$ of all subprobabilistic
distributions on $S$ a~complete lower semilattice, with the~bottom
element $0$=``no information''.

A~random variable $\alpha:S\to \BBR_+$ is called
a~\emph{probabilistic predicate}, and $\alpha(s)$ can be treated as
a~degree of appropriateness of $s\in S$ for some purpose (the~more,
the~better). In particular, if $\alpha(S)\subset\{0,1\}$, then all
elements of $S$ are divided into ``bad'' and ``good''. For
a~subprobabilistic distribution $D$, the~\emph{expectation}
$\int_D\alpha=\sum_{s\in S}D(s)\cdot\alpha(s)$ is a~maximal
expected degree guaranteed by~$D$.

A~\emph{deterministic probabilistic program} $p:S\to \bar S$ sends
each initial state $s\in S$ to a~subprobabilistic distribution
$p(s)$ of possible finite states, where the~probability
$1-\sum_{s'\in S}p(s)(s')$ is related to unknown behaviour of
the~program, in particular, to the~cases when the~program does not
terminate. Similarly, a~program $p':S\to \bar S$ refines a~program
$p:S\to\bar S$ (written $p\sqsubseteq p'$) if $p(s)\sqsubseteq
p'(s)$ for each initial state $s\in S$. If an~initial probability
distribution is partially described (estimated from below) by
a~subprobabilistic distribution $D\in \bar S$, then a~probability
of a~final state $s'\in S$ is greater or equal than
$D'(s')=\sum_{s\in S}D(s)\cdot p(s)(s')$. Therefore, for
a~probabilistic predicate $\beta:S\to\BBR_+$, the~expectation after
execution of the~program has the~best estimate from below:
$$
\int_{D'}\beta=\sum_{s,s'\in S}D(s)\cdot
p(s)(s')\cdot\beta(s').
$$
A~predicate $\alpha:S\to\BBR_+$ is called a~\emph{(probabilistic)
precondition} for $\beta$, and $\beta$ then is
a~\emph{(probabilistic) postcondition} for $\alpha$, if for
\emph{each} initial subprobabilistic distribution $D\in\bar S$ and
the~respective final subprobabilistic distribution $D'\in\bar S$,
we have $\int_D\alpha\le \int_{D'}\beta$, i.e., the~expected value
$\eps\ge 0$ of $\alpha$ guarantees that the~expectation of $\beta$
is also equal or greater than~$\eps$. It is easy to see that
the~\emph{strongest} (i.e., the~least) \emph{postcondition}
$sp(p)(\alpha)$ of $\alpha$ is determined with the~formula
$$
sp(p)(\alpha)(s')=\sum_{s\in S}\alpha(s)\cdot p(s)(s'),
\quad s'\in S.
$$
Observe that all probabilistic predicates on $S$ form a~cone, and
the~mapping $sp(p)$ is additive and positively uniform, i.e.,
preserves multiplication by non-negative numbers. In this paper we
shall construct and investigate an~analogue of this mapping.
Similarly, for a~given predicate $\beta\in \bar S$,
a~\emph{weakest} (greatest) \emph{precondition} $wp(p)(\beta)$ can
be found. See~\cite{MMS:ProbPredTrans:96} on how nondeterminism can
be incorporated into this model by mapping each initial state not
to a~single distribution, but to a~set of distributions.

This is also closely related to the~notion of
\emph{approximate correctness} of a~computer
program~\cite{MinY:ReasProbProgr:03}. Although a~number that
expresses ``approximateness'' can be also treated as degree of
belief, the~entire theory by Mingsheng Ying is based on
probabilistic logic and well suited to study probabilistic
programs. It is also focused more on uncertainty of assumptions and
conclusions than on imprecision in description of input and output
data, as one could expect based on the~term ``approximate''. For
example, the~refinement index of two probabilistic predicates is
defined as the~belief probability to which one probabilistic
predicate is refined by another. There are several parallels
between this theory and what we are doing in the~sequel.

This approach, however, has intrinsic restrictions: we assume that
a~system is sufficiently described with knowledge which states or
random events (sets of states) are realized, or what are
the~probabilities of their realization. For a~simple program, like
the~examples in~\cite{MMS:ProbPredTrans:96}, this assumption is
realistic, but if, e.g., our program removes artifacts from
a~sufficiently large colour image, then the~state space $S$ is too
huge to apply the~above apparatus. To reduce $S$, one can divide
all possible images into a~reasonable number of classes. Boundaries
between these classes cannot be clear; therefore the~predicates
will not be tolerant to small changes in images. Next, careful
study of probability distributions of the~class of a~possible
output for a~given class of an input image is a~nontrivial task.
Even if this goal is achieved, the~respective predicate
transformers describe
\emph{average} results, and say nothing about rare
extreme cases, which may make the~program unusable.

For such ``huge-dimensional'' cases we suggest to resign from
the~purely probabilistic approach and to decrease
the~``dimensionality'' by allowing fuzzy predicates. The~idea is to
have less predicates, which may be ``more or fewer'' true, and their
values for each possible portion of information about a~system
present the~greatest known degrees of truth, certainty, precision,
quality etc, which we can reliably count for. For example, such
a~predicate can assign to each square part, with integer
coordinates of the~vertices, of a~given image a~numerical measure
of its quality. Then an~image is incompletely but efficiently
described with a~finite collection of numbers, which is considered
to be the~value of the~predicate. Observe that two such collections
can be incomparable, e.g., if two images are damaged in different
places. Hence the~considered predicates can attain values in sets
which are only partially ordered, although fuzziness is most often
expressed on a~numeric scale, e.g., $[0,1]$.

>From now on we shall talk about ``truth values'' of fuzzy
predicates, but this term is used for the sake of convenience and
does not restricts possible interpretations to fuzzy logic only,
although it is also possible. We expect that all known semantics of
fuzziness
\cite{Berg:IntrManValFuzLogic:08,DubPr:GradUncBip-MakSense:10} can
be applied; see the~examples in the~next section.

Fuzzy predicate transformers also have been studied mostly in
$[0,1]$-settings~\cite{ChWu:ImplBasFuzzPredTr:05,ChWu:SemWlpSlp:11}.
This paper is devoted to constructing and investigating $L$-fuzzy
(where $L$ is a~suitable lattice) strongest postcondition predicate
transformers that are determined by state transformers, i.e., by
$L$-fuzzy knowlegde about what we can expect (more precisely, what
is guaranteed in the~worst case) for each initial state of
a~system. We are interested in order and topological properties of
predicate transformers. It will be shown that spaces of predicates
are idempotent semimodules, which are analogues of vector spaces,
and under certain (not very restrictive) conditions the~strongest
postcondition predicate transformers are linear or affine
continuous mappings between these semimodules.


\section{Semimodules of monotonic predicates}

Throughout this paper, if $f,g$ are functions with a~common domain,
$\alpha$ is a~constant, and $*$ is a~binary operation, then we
denote by $f*g$, $\alpha*f$ and $f*\alpha$ the functions with
the~same domain obtained by pointwise application of the~operation
$*$ (provided it is defined for the~corresponding values). In
the~sequel $\psup$ and $\pinf$ for a~family of functions with
a~common domain to a~poset will denote the~pointwise suprema and
infima, respectively.

See \cite{GHKLM:ContLattDom:03} for basic definitions and facts on
partially ordered sets, including continuous semilattices and
lattices. Here we shall recall only notation and a~few definitions.
For a~poset $X$, the~same set, but with the~reversed order, is
denoted by  by~$X^{op}$. An~element $a$ \emph{approximates} $b$ or is
\emph{way below} $b$, in a~poset $X$, which is written as $a\ll b$,
if, for each directed subset $C\subset X$ such that $b\le\sup C$,
there is $c\in C$ such that $a\le c$. A~poset $X$ is called
\emph{continuous} if, for each $b\in X$, the set of all $a\ll b$ is
directed and has $b$ as its lowest upper bound. A~poset is
\emph{directed complete} if each its non-empty \emph{directed}
subset has a~least upper bound. A~continuous directed complete
poset is called a~\emph{domain}. A~domain which is additionally
a~meet semilattice (a~complete lattice) is called
a~\emph{continuous semilattice} (respecticely a~\emph{continuous
lattice}).

The~\emph{Scott topology} on a~poset $X$ is the~least topology such
that all lower sets $C$ that are closed under directed suprema are
closed. The~\emph{lower topology} on $X$ is the~least topology such
that the~sets $\{a\in X\mid b\le a\}$ are closed for all $b\in X$.
The~join, i.e., the~least topology that contains the~Scott and
the~lower topologies, is called the~\emph{Lawson topology}.

In the sequel $L$ will be a~completely distributive lattice, i.e.,
a~compact Hausdorff distributive Lawson lattice with its Lawson
topology. A~topological lattice (semilattice) is said to be Lawson
if for each point it possesses a~local base that consists of
sublattices (respectively of subsemilattices). Note that the~same
is true for $L^{op}$. We denote by $0$, $1$, $\oplus$, and
$\otimes$ the~bottom element, the~top element, the~join, and
the~meet in $L$, respectively. The elements of this (arbitrary, but
fixed throughout the~paper) lattice will be used to express truth
values. The~operation $\oplus$ is the~disjunction, but
the~conjuction does not necessarily coincide with~$\otimes$.
Although complete distributivity is a~very strong requirement,
a~lot of important lattices fall into this class, e.g., all
complete linearly ordered sets, including $I=[0,1]$ or any other
segment in $\BBR$, all finite distributive lattices, all products
of completely distributive lattices, in particular, $I^\tau$ for
all cardinals $\tau$. In fact, a~lattice is completely distributive
if and only if it is order isomorphic to a~complete sublattice of
some $I^\tau$.

We shall also use basic notions of denotational semantics of
programming languages. Consider a~state of a~computational process
or a~system. All possible (probably incomplete) portions of
information we can have about this state form a~\emph{domain of
computation}~$D$~\cite{Eda:DomComput:97}. This set carries
a~partial order $\le$ which represents a~hierarchy of information
or knowledge: the~more information an~element contains (i.e.,
the~more specific/restrictive it is), the~higher it is.
See~\cite{Eda:DomComput:97} for more details, in particular, for
an~explanation why it is natural to demand that $D$ is a~domain, i.e.,
a~continuous directed complete poset. In addition to this, it is
also often required that there is a~least element $0\in D$ (no
information at all), and that for all $a$ and $b$ in $D$ there is
a~meet $a\land b$, which, e.g., can be (but not necessarily is)
treated as ``$a$ or $b$ is true''.

Following \cite{HH:DualTheoQuantSem:98}, for a~domain $D$ we call
elements of the~set $[D\to L^{op}]^{op}$ $L$-\emph{fuzzy monotonic
predicates} on $D$ (here $[A\to B]$ stands for the~set of mappings
from $A$ to $B$ that are Scott continuous, i.e., they preserve
directed suprema). For $m\in [D\to L^{op}]^{op}$ and $a\in D$, we
regard $m(a)$ as the~truth value of $a$; hence it is required that
$m(b)\le m(a)$ for all $a\le b$. The~second ${}^{op}$ means that we
order fuzzy predicates pointwise, i.e., $m_1\le m_2$ iff $m_1(a)\le
m_2(a)$ in $L$ (not in $L^{op}$~!) for all $a\in D$. We denote
$\uM_{[L]}D=[D\to L^{op}]^{op}$, and, for a~domain $D$ with
a~bottom element, consider also the~subset $M_{[L]}D\subset
\uM_{[L]}D$ of all
\emph{normalized} predicates that take $0\in D$ (no information) to
$1\in L$ (complete truth). Observe that $M_{[L]}D$ is a~complete
sublattice of~$\uM_{[L]}D$.

\begin{example}
Let a~system have a~finite or countable state space $S$. Each
subset $A\subset S$ is identified with it characteristic mapping
$\chi_A:S\to\{0,1\}$, which is a~Boolean predicate ``current state
$s$ is in~$A$''. A~smaller subset $A$ corresponds to more
information; therefore the~set $D$ of all subsets of $S$ is
partially ordered by reverse inclusion. Then $D$ is a~continuous
lattice, and the~$\{0,1\}$-fuzzy monotonic predicates on $D$ are
precisely $\chi_A$ for all~$A\subset S$.

If the~system changes its state randomly, then different schemes
are possible. Generally, an~incomplete probabilistic knowledge is
a~mapping $m:D\to [0,1]$ such that for all $A\subset S$
the~probability $P(A)$ is at least $m(A)$. Of course, $A\le B$,
i.e., $A\supset B$, implies $m(A)\ge m(B)$, and $\sigma$-additivity
of probability requires that $m$ sends the~directed unions
of subsets of~$S$ to the~corresponding suprema in~$[0,1]$. Thus $m$
is a~$[0,1]$-fuzzy monotonic predicate.

Observe that $m$ may not necessarily be reduced to a~collection of
estimates for the~probabilities of individual states $s\in S$. For
example, if all that we know is $P(\{s_1,s_2\})\ge 0{,}5$, then
the~only subprobabilistic distribution that is surely less or equal
than the~actual distribution is trivial, i.e., zero for all states.

Of course, $m$ can be determined by (sub)probabilistic
distributions. Let an~exact probability distribution be unknown,
but one of $n$ possible, which are bounded from below respectively
by subprobabilistic distributions $P_1,P_2,\dots,P_n\in\bar S$.
The~greatest guaranteed probability of a~random event $A\in D$ is
equal to $m(A)=\inf_{1\le i\le n}\sum_{s\in A}P_i(s)$. Then $m$ is
a~$[0,1]$-valued fuzzy monotonic predicate, which ``aggregates''
all possible probability distributions in the~assumption of
``demonic'' non-determinism.

Thus numeric fuzzy predicates can arise in purely probabilistic
settings, with the~semantics ``truth value = guaranteed
probability''. Observe that the~probability of $S$ is always $1$,
hence the~mentioned predicates may be considered normalized.
\end{example}

\begin{example}\label{ex.m-scale-pred}
Let an~image be divided into $n$ parts, and the~quality of each of
them can be rated in the~scale $L=\{0,1,\dots,m\}$, e.g.,
$0$=``awful'', $1$=``bad'', \dots, $m$=``perfect''. Then the~state
space is equal to $S=L^n$. The~domain of computation $D$ can also
be put equal to $L^n$, and $d=(d_1,d_2,\dots,d_n)$ will mean
``the~actual quality $s_i$ of $i$-th part is not worse than $d_i$
for all $1\le i\le n$''. This implies that $(d_1,d_2,\dots,d_n)\le
(d'_1,d'_2,\dots,d'_n)$ in $D$ if and only if $d_1\le d_1'$,
$d_2\le d_2'$, \dots, $d_n\le d_n'$.

For each $q=(q_1,q_2,\dots,q_n)\in L^n$, let the~predicates
$m_q,m'_q,m''_q:D\to L$ be defined by the~formulae:
\begin{gather*}
m_q\bigl((d_1,d_2,\dots,d_n)\bigr)=
\max\bigl\{k\in L\mid d_i\ge q_i-(m-k)
\text{ for all }i=1,2,\dots,n\bigr\},
\\
m'_q\bigl((d_1,d_2,\dots,d_n)\bigr)=
\max\bigl\{k\in L\mid d_i\ge\min\{k,q_i\}
\text{ for all }i=1,2,\dots,n\bigr\},
\\
m''_q\bigl((d_1,d_2,\dots,d_n)\bigr)=
\max\bigl\{k\in L\mid \max\{d_i,m-k\}\ge q_i
\text{ for all }i=1,2,\dots,n\bigr\}
\end{gather*}
for all $(d_1,d_2,\dots,d_n)\in S$. Then
$m_q\bigl((d_1,d_2,\dots,d_n)\bigr)$ shows the~worse relative loss
of quality w.r.t.\ $(q_1,q_2,\dots,q_n)$,
$m'_q\bigl((d_1,d_2,\dots,d_n)\bigr)$ shows ``below what degree''
the~quality of $(d_1,d_2,\dots,d_n)$ is not worse than
$(q_1,q_2,\dots,q_n)$, and $m''_q\bigl((d_1,d_2,\dots,d_n)\bigr)$
shows ``above what degree'' the~quality of $(d_1,d_2,\dots,d_n)$ is
not worse than $(q_1,q_2,\dots,q_n)$. In all these cases
the~predicates compare the~guaranteed quality of an~input with
a~desired one. Thus we can construct a~predicate like ``the~image
is perfect at the~center and at least good at the~angles''.

Moreover, we can rate parts of an~image in several aspects, with
separate scales $L_1,L_2,\dots,L_r$ for each, then the~resulting
$L=L_1\times L_2\times\dots\times L_r$ will be a~finite
distributive lattice, which is not linearly ordered.
\end{example}

\medskip

It follows from~\cite[Theorem~4]{Erne-ZdistrFunSp:98} (classified
as ``folklore knowledge'' in \cite{HH:DualTheoQuantSem:98}) that,
for a~domain $D$ and a~completely distributive lattice $L$, the~set
$[D\to L^{op}]$ is a~completely distributive lattice as well. Hence
this is also valid for $\uM_{[L]}D$ and (if $D$ contains a~least
element) $M_{[L]}D$.

For an element $d_0\in D$, we denote by $\eta_{[L]} D(d_0)$
the~function $D\to L$ that sends each $d\in D$ to $1$ if $d\le d_0$
and to $0$ otherwise. It is easy to see that $\eta_{[L]} D(d_0)\in
M_{[L]}D\subset \uM_{[L]}D$, and $\delta^D_L=\eta_{[L]} D(0)$ is
a~least element of $M_{[L]}D$.

\begin{lemma}\label{lem.eta-cont}
For a~domain $D$, the mapping $\eta_{[L]} D:D\to \uM_{[L]}D$ is
continuous w.r.t.\ the~Scott topologies and w.r.t.\ the~lower
topologies. If $D$ is a~complete continuous semilattice, then
$\eta_{[L]} D$ is an~embedding w.r.t.\ the~Scott topologies,
the~lower topologies, and the~Lawson topologies.
\end{lemma}

\begin{proof}
Obviously, $\eta_{[L]}D(d_1)\le \eta_{[L]}D(d_2)$ if and only if
$d_1\le d_2$. Observe also that $\eta_{[L]} D(d_0)$ is a~least
$m\in\uM_{[L]}D$ such that $m(d_0)=1$. If $\CCD\subset D$ is
directed and $\sup\CCD=d_0$, then $\sup\{\eta_{[L]}D(d)\mid
d\in\CCD\}$ is a~least $m\in\uM_{[L]}D$ such that $m\ge
\eta_{[L]}D(d)$ for all $d\in\CCD$, which is equivalent to $m(d)=1$
for all $d\in\CCD$. Since $m:D\to L^{op}$ is Scott continuous,
i.e., it preserves directed suprema, which is, in turn, equivalent
to $m(\sup\CCD)=m(d_0)=1$. By the~above such $m$ is equal to
$\eta_{[L]}D(d_0)$. Hence $\eta_{[L]}D$ preserves directed suprema
as well.

To show that $\eta_{[L]}D$ is lower continuous, it suffices to show
that, for all $m\in \uM_{[L]}D$, the set
$$
\eta_{[L]}D{}^{-1}(\{m\}\ups)=
\{d_0\in D\mid \eta_{[L]}D(d_0)\ge m\}
$$
is closed in the lower topology on $D$. The~inequality
$\eta_{[L]}D(d_0)\ge m$ means that $\eta_{[L]}D(d_0)(d)=1$ for all
$d\in D$ such that $m(d)\ne 0$; in other words, $d_0$ is an~upper
bound of the set $\{d\in D\mid m(d)\ne 0\}$. This implies that
$$
\eta_{[L]}D{}^{-1}(\{m\}\ups)=
\bigcap\bigl\{\{d\}\ups\subset D\mid m(d)\ne 0, d\in D\bigr\},
$$
which is closed in the~lower topology on $D$.

If $D$ is a~complete continuous semilattice, then it is compact
Hausdorff in its Lawson topology; therefore a~continuous injective
mapping from it to a~compactum $\uM_{[L]}D$ is an~embedding. Due to
the~completeness of $D$, this implies that the~isotone mapping
$\eta_{[L]}D$ is also an~embedding w.r.t.\ the~Scott topologies and
w.r.t.\ the~lower topologies.
\end{proof}

Therefore we consider $D$ as a~sub\emph{dcpo} of $\uM_{[L]}D$, and
a~complete continuous semilattice $D$ is additionally a~subspace of
$\uM_{[L]}D$ w.r.t.\ the~Scott, the~lower, and the~Lawson
topologies on the~both sets.

Infima and finite suprema in the~complete lattices $\uM_{[L]}D$ and
$M_{[L]}D$ of functions are taken pointwise, whereas arbitrary
suprema are described by the~following easy, but useful statement.
For a~function $f:D\to L$, let
$$
f^u(b)=\inf\{f(a)\mid a\in D,a\ll b\},\text{ for all }b\in D.
$$
Observe that $f^u$ is always a~monotonic predicate. Moreover
\cite[Lemma~I.4]{Verw:RandUSC:97}:

\begin{lemma}
For an~antitone function $f:D\to L$, the~function $f^u$ is
the~least monotonic predicate $f'$ such that $f\le f'$ pointwise.
\end{lemma}

Hence, for a~family $\CCF\subset \uM_{[L]}D$ (or $\CCF\subset
M_{[L]}D$), we have $\inf\CCF=\pinf\CCF$, $\sup\CCF=(\psup\CCF)^u$.
For finite $\CCF$, the~latter ${}^u$ can be dropped.

\begin{lemma}\label{lem.comp-psup}
Let a set $\CCF\subset \uM_{[L]}D$ (or $\CCF\subset M_{[L]}D$) be
compact in the~relative lower topology. Then $\psup\CCF\in
\uM_{[L]}D$ (resp.\ $\psup\CCF\in M_{[L]}D$); therefore
$\sup\CCF=\psup\CCF$.
\end{lemma}

\begin{proof}
Assume to the~contrary, that there exists $a_0\in D$ such that
$$
\sup\{f(a)\mid f\in\CCF\}\ge\alpha\not\le\alpha_0=\sup\{f(a_0)\mid
f\in\CCF\}
$$
for all $a\in D$, $a\ll a_0$. The complete distributivity of $L$
implies that there is $\beta\in L$ such that $\beta\le \alpha$,
$\beta\not\le \alpha_0$, and if $\Gamma\subset L$ satisfies $\sup
\Gamma\ge\alpha$, then there is $\gamma\in \Gamma$, $\gamma\ge\beta$
(such $\beta$ is said to be \emph{way-way below} $\alpha$,
cf.~\cite{GHKLM:ContLattDom:03}). The~set
$$
\CCF_a=
\{f\in \CCF\mid f(a)\ge \beta\}=
\{f\in \CCF\mid f\ge m\},
$$
where
$$
m(a')=
\begin{cases}
\beta, a'\le a,\\
0, a'\not\le a,
\end{cases}
\text{ for }a'\in D,
$$
is closed in $\CCF$. The~family $\{\CCF_a\mid a\ll a_0\}$ of
nonempty sets is directed; therefore by compactness it has a~common
element $f_0\in\CCF$, i.e., $f_0(a)\ge \beta$ for all $a\ll a_0$.
Then by the Scott continuity of $f_0:D\to L^{op}$ we obtain
$$
\alpha_0=\sup\{f(a_0)\mid f\in\CCF\}\ge f_0(a_0)\ge\beta,
$$
which is a~contradiction.
\end{proof}

We use notation $\bp$ and $\bm$ for respectively joins and meets
both in $\uM_{[L]}D$ and $M_{[L]}D$.

In the sequel we shall additionally require that $L$ be
a~\emph{unital quantale}~\cite{Rosen:Quantales:90}, i.e., there
exists an~associative binary operation $*:L\times L\to L$ such that
$1$ is a~two-sided unit and $*$ is infinitely distributive w.r.t.\
supremum in both variables, which is equivalent to being continuous
w.r.t.\ the~Scott topology on $L$. Observe that, for such $*$, its
infinite distributivity w.r.t.\ infima also means continuity
w.r.t.\ the~Lawson topology on $L$. Recall that we treat $\op$ as
a~disjunction, and $*$ will be a~(possibly noncommutative)
conjunction in an $L$-valued fuzzy
logic~\cite{Haj:FuzNoncomLog:03}. The~Boolean case is obtained for
$L=\{0,1\}$, $\op=\lor$ and $*=\land$. On the~other hand, let
the~finite linearly ordered set $L=\{0,1,\dots,m\}$ be used to
express absolute and relative quality of input, certainty,
precision, etc., cf.\ Example~\ref{ex.m-scale-pred}. Then
the~operations $i*j\equiv \min\{i,j\}$ and
$i*j\equiv\max\{i+j-m,0\}$ can be reasonable choices, which reflect
the~natural assumption that combination of two distorted,
imprecise, or uncertain inputs produces an~equally or more
distorted, imprecise, or uncertain output.

\begin{lemma}\label{lem.alpha-u}
For $\alpha\in L$, a predicate $m\in \uM_{[L]}D$, and an~antitone
function $f:D\to L$, we have $m(b)\ge \alpha*f(b)$ (resp.\ $m(b)\ge
f(b)*\alpha$) for all $b\in D$ if and only if $m(b)\ge
\alpha*f^u(b)$ (resp.\ $m(b)\ge f^u(b)*\alpha$) for all $b\in
D$.
\end{lemma}

\begin{proof}
Since $f\le f^u$, ``if'' is trivial. Assume that $m(b)\ge
\alpha*f(b)$ for all $b\in D$. Then for all $a\in D$, $a\ll b$
the~inequality $f^u(a)\ge f(b)$ implies $m(a)\ge \alpha*f^u(b)$.
Putting $a\to b$, we obtain $m(b)\ge \alpha*f^u(b)$.
\end{proof}

\begin{remark}\label{rem.alpha-f-u}
The~latter statement can be expressed by the~formulae:
$$
(\alpha*f)^u=(\alpha*f^u)^u,\quad (f*\alpha)^u=(f^u*\alpha)^u,
$$
for each antitone function $f:D\to L$ and $\alpha\in L$. It is also
easy to see that, for a~family $\{f_i\mid i\in\CCI\}$ of antitone
functions $D\to L$, the~equality
$$
\bigl(\psup_{i\in\CCI}f_i\bigr)^u = \bigl(\psup_{i\in\CCI}(f_i)^u\bigr)^u
$$
is valid.
\end{remark}

The operation $*$ induces binary operations $\bd$ and $\bs$ on
the~posets $\uM_{[L]}D$ and $M_{[L]}D$, which make them
$L$-\emph{idempotent compact Lawson
semimodules}~\cite{Nyk:ContDualContLsem:12}. Recall that a (left
idempotent) $(L,\op,*)$-semimodule~\cite{Akian:DensInvMeas:99} is a
set $X$ with operations $\bp:X\times X\to X$ and $\os:L\times X\to
X$ such that for all $x,y,z\in X$, $\alpha,\beta\in L$~:

(1) $x\bp y=y\bp x$;

(2) $(x\bp y)\bp z=x\bp (y\bp z)$;

(3) there is an (obviously unique) element $\bar 0\in X$ such that
$x\bp\bar 0=x$ for all $x$;

(4) $\alpha\os(x\bp y)=(\alpha\os x)\bp(\alpha\os y)$,
$(\alpha\op\beta)\os x=(\alpha\os x)\bp(\beta\os x)$;

(5) $(\alpha*\beta)\os x=\alpha\os(\beta\os x)$;

(6) $1\os x=x$; and

(7) $0\os x=\bar 0$.

Observe that these axioms imply that $(X,\bp)$ is an upper
semilattice with a bottom element $\bar 0$, and $\alpha\os\bar
0=\bar 0$ for all $\alpha\in L$. The~operation $\os$ is isotone in
both variables.

Hence an~$(L,\op,*)$-semimodule is an~analogue of a~vector space.
Similarly, analogues exist for linear and affine mappings.
A~mapping $f:X\to Y$ between $(L,\op,*)$-semimodules is called
\emph{linear} if, for all $x_1,\dots,x_n\in X$ and
$\alpha_1,\dots,\alpha_n\in L$ , the equality
$$
f(\alpha_1\os x_1\bp\dots\bp\alpha_n\os x_n)=
\alpha_1\os f(x_1)\bp\dots\bp\alpha_n\os f(x_n)
$$
is valid. If the latter equality is ensured only whenever
$\alpha_1\op\dots\op \alpha_n=1$, then $f$ is called
\emph{affine}. Observe that an~affine mapping $f$ preserves joins,
i.e., $f(x_1\bp x_2)=f(x_1)\bp f(x_2)$ for all $x_1,x_2\in X$.
An~affine mapping is linear if and only if it preserves the~least
element.

We call a~triple $(X,\bp,\os)$ a~\emph{continuous
$(L,\op,*)$-semimodule}~\cite{Nyk:ContDualContLsem:12} if
$(X,\bp,\os)$ is an~$(L,\op,*)$-semimodule, $X$ is a~continuous
(hence complete) lattice, and $\os:L\times X\to X$ is infinitely
distributive w.r.t.\ all suprema in both variables. Then $X$ with
its Lawson topology is a~compact Hausdorff Lawson lower semilattice
with a~top element, and $\os$ is jointly continuous w.r.t.\
the~Scott topologies on $L$ and $X$.

For $m\in \uM_{[L]}D$, we define $\alpha\bd m$ to be a~least
predicate $m':D\to L$ such that $\alpha*m(b)\le m'(b)$ for all
$b\in D$, i.e., $\alpha\bd m=(\alpha*m)^u$. Then:
$$
(\alpha\bd m)(b)=\inf\{\alpha*m(a)\mid a\in D,a\ll b\}.
$$
For $m\in M_{[L]}D$, we need to ``adjust'' the~result:
$$
(\alpha\bs m)(b)= (\alpha\bd m)(d)\bp \delta^D_L=
\begin{cases}
(\alpha\bd m)(b),&b\ne 0;\\
1, &b=0.
\end{cases}
$$

\begin{lemma}
For $\alpha,\beta\in L$, $m\in \uM_{[L]}D$:
$$
\alpha\bd(\beta\bd m)=(\alpha*\beta)\bd m.
$$
\end{lemma}

\begin{proof}
By Remark~\ref{rem.alpha-f-u}:
$$
\alpha\bd(\beta\bd m)=\bigl(\alpha*(\beta*m)^u\bigr)^u=
\bigl(\alpha*(\beta*m)\bigr)^u=(\alpha*\beta)\bd m.
$$
\end{proof}

Now the~equality
$$
\alpha\bs(\beta\bs m)=(\alpha*\beta)\bs m
$$
for all $\alpha,\beta\in L$, $m\in M_{[L]}D$ is immediate. Both
operations $\bd$ and $\bs$ are infinitely distributive w.r.t.\
supremum in the~both arguments (because $*$ is such an~operation);
hence, both are lower semicontinuous. Using routine, but
straightforward calculations (\cite{Nyk:ContDualContLsem:12}; the
same but in terms of hyperspaces in \cite{Nyk:AdMonCompLawsSem:11})
we obtain:

\begin{proposition}
The~triples $(\uM_{[L]}D,\bp,\bd)$ and $(M_{[L]}D,\bp,\bs)$ are
continuous $(L,\op,*)$-semimodules.
\end{proposition}

\begin{remark}
It is easy to see that, if $*$ is also infinitely distributive
w.r.t.\ infimum, then $\alpha*m\in [D\to L^{op}]^{op}$ for all
$\alpha\in L$, $m\in [D\to L^{op}]^{op}$. Therefore, in this case
$\alpha\bd m$ coincides with $\alpha*m$.
\end{remark}

For two~predicates $m_1,m_2:D\to L$, their join (i.e.,
the~argumentwise supremum) $m_1\bp m_2$ can be interpreted as
disjunction: ``$m_1$ or $m_2$''. Multiplication of a~predicate
$m:D\to L$ by $\alpha\in L$ either does not change this predicate
or makes it more ``pessimistic'', or, equivalently, more
``demanding''. Since the~sets of $L$-fuzzy monotonic predicates are
``vector-like'' spaces, we can apply to them the~tools of
idempotent linear algebra and idempotent functional analysis,
although these theories are rather limited and poor comparing to
the~``conventional'' classical analogues. In particular, results of
\cite{Nyk:ContDualContLsem:12} allow:
\begin{itemize}
\item
to approximate $L$-fuzzy
monotonic predicates from below and from above with predicates that
attain only finite sets of values;
\item
to study and approximate
predicates with special properties, e.g., meet- and
join-preserving; and
\item
to construct the~predicate that is dual to a~given
one, if the latter expresses an~undesirable property which have to
be avoided, etc.
\end{itemize}

\section{Strongest postcondition predicate transformers}

We treat each~mapping $m:D\to L$ as ``it is known that, for each
$d\in D$, its truth value is at least $m(d)$''. Similarly,
an~arbitrary mapping $\phi:D\to \uM_{[L]}D'$ is interpreted as ``if
$a\in D$ is true, then the~truth value of each $b\in D'$ is at
least $\phi(a)(b)$''. Note that $\phi(a)(b)$ is implicitly
considered as a~``conditional'' truth value, i.e., if $a$ is
``partially true'' at a~degree $\ge\alpha$, then $b$ is true at
least at a~degree $\alpha*\phi(a)(b)$.

Hence, such a~$\phi$ is an~$L$-fuzzy \emph{state transformer}. For
a~given $\phi$, we say that $m:D\to L$ is a~\emph{precondition} and
$m':D'\to L$ is a~\emph{postcondition} for each other w.r.t.\
$\phi$, if, for all $a\in D$ and $b\in D'$, the~``guaranteed''
truth value $m'(b)$ is greater or equal to $m(a)*\phi(a)(b)$, i.e.,
to the~result of \emph{modus ponens}.

Obviously, for an~antitone function $m:D\to L$, its
\emph{strongest} (\emph{least}) postcondition $\usp(\phi)(m)$ in
$\uM_{[L]}D'$ is determined by the~equality
$$
\usp(\phi)(m)(b)=\inf\bigl\{\sup\{m(a)*\phi(a)(b') \mid a\in D\}
\mid b'\in D',b'\ll b\bigr\}, b\in D'.
$$

Again, if we restrict ourselves to normalized predicates,
the~strongest postcondition must be corrected:
$$
sp(\phi)(m)(b)=
\usp(\phi)(m)(b)\bp \delta^D_L=
\begin{cases}
\usp(\phi)(m)(b),&b\ne 0;\\
1, &b=0.
\end{cases}
$$

It is easy to see that, for all $d\in D$ and isotone
$\phi:D\to\uM_{[L]}D'$, we have
$\usp(\phi)(\eta_{[L]}D(d))=\phi(d)$, hence $\usp(\phi)$ is
an~isotone extension of $\phi$. Similarly, for an~isotone mapping
$\phi:D\to M_{[L]}D'$, the~mapping $sp(\phi)(\eta_{[L]}D(d))$ is
an~isotone extension as well. The mapping $\usp(\phi)$ and
$sp(\phi)$ are called ($L$-fuzzy)
\emph{strongest postcondition predicate transformers} induced by
the~state transformer $\phi$, and are analogues of crisp (i.e.,
Boolean) predicate transformers, which were introduced by
Dijkstra~\cite{Dijk:Guarded:75}. Compare also with
the~\emph{weakest precondition predicate transformers},
cf.~\cite{ChWu:ImplBasFuzzPredTr:05,ChWu:SemWlpSlp:11}. Their
$L$-valued ``angelic'' and ``demonic'' analogues were introduced
and investigated in~\cite{DMR:LattValPredTrIntSys:10} by means of
topology. The~latter reference contains also an~example of
a~security system, which analyzes security threats of different
severity and nature and imposes security measures of
the~corresponding level. This is naturally expressed with elements
of lattices; therefore the~authors propose to ``consider possible
definitions for lattice-valued predicate transformers''. Here is
another example.

\begin{example}
Assume that a~program processes a~sequence of $n$ frames. The~quality $s_i$
of $i$-th frame is rated in the~scale $L=\{0,1,\dots,m\}$.
The~domain of computation is equal to $D=L^n$, and the~meaning of
$d=(d_1,d_2,\dots,d_m)$ is ``$s_1\ge d_1$, $s_2\ge d_2$, \dots,
$s_n\ge d_n$''. The~multiplication $i*j=\max\{i+j-m,0\}$ is
considered on $L$, making it a~finite quantale. The~truth value of
$d=(d_1,d_2,\dots,d_n)$ is defined as
$$
\max\{k\in L\mid s_i\ge d_i*k\text{ for all }i=1,2,\dots,n\}
$$
(observe that it is $m_s(d)$ for $s=(s_1,s_2,\dots,s_n)$, cf.\
Example~\ref{ex.m-scale-pred}). Assume that it is known that, if
the~quality of $i$-th frame, $0<i<n$, is $\ge k-1$, and the~quality
of the~two neighboring frames is $\ge k$, then, after the~program
execution, the~quality of $i$-th frame will be $\ge k$, for all
$1\le k\le m$. This information can be expressed via the~state
transformer $\phi:D\to \uM_{[L]}D$ that sends
$$
s=(0,\dots,
\underset{i-1}{{m}},
\underset{i}{{m-1}},
\underset{i+1}{{m}},
\dots,0),\text{ for }0<i<n,
$$
to $m_q$, where
$$
q=(0,\dots,
\underset{i-1}{{0}},
\underset{i}{{m}},
\underset{i+1}{{0}},
\dots,0),
$$
and all other $s\in D$ to the~constant zero predicate. Similarly we
can add the~fact that the~quality of each frame will not be worse
than before, etc. The~resulting predicate transformer $\usp(\phi):
\uM_{[L]}D\to \uM_{[L]}D$ sends a~known quality of the~frames
\emph{before} the~program run to the~most guaranteed quality
\emph{after} its execution.
\end{example}

To simplify our exposition, we consider in this section not
necessarily normalized monotonic predicates.

\begin{lemma}\label{lem.post-u}
With respect to a~Scott continuous mapping $\phi:D\to\uM_{[L]}D'$,
a~monotonic predicate $m':D'\to L$ is a~postcondition for
an~antitone function $m:D\to L$ if and only if $m'$ is
a~postcondition for $m^u:D\to L$.
\end{lemma}

\begin{proof}
Since $m\le m^u$, ``if'' is immediate. Let $m'(b)\ge
m(a)*\phi(a)(b)$ for all $a\in D$, $b\in D'$. Then $m'(b)\ge
m(a')*\phi(a')(b)\ge m^u(a)*\phi(a')(b)$ for all $a'\ll a$. This
implies $m'\ge m^u(a)*\psup_{a'\ll a}\phi(a')$, therefore by
Lemma~\ref{lem.alpha-u}
$$
m'\ge m^u(a)*(\psup_{a'\ll a}\phi(a'))^u=m^u(a)*\sup_{a'\ll
a}\phi(a')=m^u(a)*\phi(a'),
$$
the~last equality is due to the~Scott continuity of $\phi$.
\end{proof}

\begin{proposition}\label{stat.scott-sp}
Let $\phi$ be a~mapping $D\to\uM_{[L]}D'$. Then $\usp(\phi):
\uM_{[L]}D\to \uM_{[L]}D'$ preserves joins (hence finite suprema).
For an~isotone $\phi$, the~mapping $\usp(\phi)$ preserves all
suprema if and only if $\phi$ is Scott continuous, i.e., preserves
directed suprema.
\end{proposition}

\begin{proof}
Let $m=m_1\bp m_2$, for $m,m_1,m_2\in \uM_{[L]}D$. Then, for $m'\in
\uM_{[L]}D'$, $a\in D$, $b\in D'$, the~inequality
$m'(b)\ge (m_1\bp m_2)(a)*\phi(a)(b)$ is valid if and only if both
$m'(b)\ge m_1(a)*\phi(a)(b)$ and $m'(b)\ge m_2(a)*\phi(a)(b)$ are
satisfied. Therefore
\begin{gather*}
\min\{m'\in \uM_{[L]}D'\mid
m'(b)\ge (m_1\bp m_2)(a)*\phi(a)(b)
\text{ for all }a\in D, b\in D'
\}=
\\
\min\{m'\in \uM_{[L]}D'\mid
m'(b)\ge m_1(a)*\phi(a)(b)
\text{ for all }a\in D, b\in D'
\}\bp
\\
\min\{m'\in \uM_{[L]}D'\mid
m'(b)\ge m_2(a)*\phi(a)(b)
\text{ for all }a\in D, b\in D'
\},
\end{gather*}
i.e.,
$$
\usp(\phi)(m_1\bp m_2)=\usp(\phi)(m_1)\bp \usp(\phi)(m_2).
$$

Now let $\phi$ be isotone. If $\usp(\phi)$ preserves all suprema,
than it is Scott continuous, as well as $\phi=\usp(\phi)\circ
\eta_{[L]}D'$.

If $\phi$ is Scott continuous and $\{m_i\mid i\in \CCI\}\subset
\uM_{[L]}D$, then due to monotonicity
$$
\usp(\phi)(\sup\limits_{i\in\CCI}m_i)\ge
\sup\limits_{i\in\CCI}\usp(\phi)(m_i).
$$

On the~other hand, $\sup\limits_{i\in\CCI}\usp(\phi)(m_i)$ is
a~postcondition for all $m_i$; hence by Lemma~\ref{lem.post-u} for
$(\psup\limits_{i\in\CCI}m_i)^u=\sup\limits_{i\in\CCI}m_i$.
Therefore
$$
\sup\limits_{i\in\CCI}\usp(\phi)(m_i)\ge
\usp(\phi)(\sup\limits_{i\in\CCI}m_i),
$$
and $\usp(\phi)$ preserves all suprema.
\end{proof}

Unfortunately, an~analogue of Proposition~\ref{stat.scott-sp} for
lower topologies is not valid, even if $*$ is infinitely
distributive w.r.t.\ both suprema and infima.

\begin{example}
Let $D=\{0,1,1'\}\cup\{1+\frac1n\mid n=1,2,3,\dots\}$ with
the~usual numeric order, except that $1'$ is an~extra copy of $1$,
and $1$ and $1'$ are incomparable. Each directed set in $D$ has
a~greatest element, hence $D$ is a directed complete continuous
poset. Thus $D$ is an~incomplete continuous semilattice with
a~least element $0$. All upper sets in $D$ are lower closed and
Scott open; therefore all~isotone mappings from $D$ to any poset
are continuous w.r.t.\ both the~lower and the~Scott topologies.

Also, let $L=D'=\{0,1\}$; $*={\land}$; and $\phi:D\to\uM_{[L]}D'$
be an~isotone mapping defined as follows:
$$
\phi(d)=\begin{cases}
\bar 0, d\in\{0,1,1'\},\\
\delta^{D'}_L,d\notin\{0,1,1'\},
\end{cases}
d\in D.
$$

Then
$$
\usp(\phi)(m)(0)=
\begin{cases}
1\text{ if there is }d\in\{1+\frac1n\mid n=1,2,3,\dots\}, m(d)=1,
\\
0\text{ otherwise}.
\end{cases}
$$

Therefore there is a~greatest element $m_1$ in the~complement of
the~preimage $\usp(\phi)^{-1}(\{\delta^{D'}_L\}\ups)$ in
$\uM_{[L]}D$:
$$
m_1(d)=
\begin{cases}
1, d\in\{0,1,1'\},\\
0, d\notin\{0,1,1'\},
\end{cases}
d\in D.
$$
Hovewer, there are no minimal elements in the~preimage itself;
hence it is not lower closed.

Thus $\usp(\phi)(m)(0)$ is not lower continuous.
\end{example}

To obtain the~required analogue, we must apply additional
requirements.

\begin{proposition}\label{stat.low-sp}
Let $D$ and $D'$ be complete continuous semilattices, $\phi:D\to
\uM_{[L]}D'$ an~isotone mapping, and $*:L\times L\to L$ infinitely
distributive also w.r.t.\ infimum in both variables. Then
$\usp(\phi)$ is lower continuous if and only if $\phi$ is lower
continuous, and in this case $\usp(\phi)$ is defined by a~simpler
formula:
$$
\usp(\phi)(m)(b)=\sup\{m(a)*\phi(a)(b) \mid a\in D\}, b\in D'.
$$
\end{proposition}

\begin{proof}
Recall that such an~operation $*:L\times L\to L$ is continuous
w.r.t.\ the~lower and the~Lawson topologies on $L$, while
the~previously required infinite distributivity w.r.t.\ supremum
implies only the~Scott continuity of~$\phi$. The~semilattices $D$
and $D'$ with the~Lawson topologies are compact Hausdorff
topological semilattices.

{\sl Necessity} is due to Lemma~\ref{lem.eta-cont}, because
$\phi=\usp(\phi)\circ \eta_{[L]} D$, and $\eta_{[L]} D$ is lower
continuous.

{\sl Sufficiency.} The~mapping that sends each $a\in D$ to
$m(a)*\phi(a)\in \uM_{[L]}D$ is continuous w.r.t.\ the~Lawson
topology on $D$ and the~lower topology on $\uM_{[L]}D$. Hence
the~set $\{m(a)*\phi(a)\mid a\in D\}$ is compact in the~lower
topology on $\uM_{[L]}D$. By Lemma~\ref{lem.comp-psup} its
pointwise limit is in $\uM_{[L]}D$; therefore it coincides with
$\usp(\phi)(m)$.

Let $m\in \uM_{[L]}D\setminus
\usp(\phi)^{-1}(\{m'\}\ups)$, $m'\in
\uM_{[L]}D'$, then
$
\usp(\phi)(m)(b)=\sup\{m(a)*\phi(a)(b)\mid a\in D\}
=\gamma\not \ge m'(b)
$
for some $b\in D'$.

The~set $\{(m(a),\phi(a)(b))\mid a\in D\}$ is contained in
the~closed, therefore compact, lower set $\{(\alpha,\beta)\in
L\times L\mid\alpha*\beta\le\gamma\}$. The~operation $*$ is isotone
and Lawson continuous. Hence there are
$\alpha_1,\beta_1,\dots,\alpha_n,\beta_n\in L$ such that the~open
set
$$
U=(L\times L)\setminus (
\{\alpha_1\}\ups\times \{\beta_1\}\ups
\cup
\dots
\cup
\{\alpha_n\}\ups\times \{\beta_n\}\ups)
$$
contains
$$
\{(\alpha,\beta)\in L\times L\mid
\alpha*\beta\le\gamma\},
$$
and $\sup\{\alpha*\beta\mid (\alpha,\beta)\in U\}=\gamma'\not
\ge m'(b)$. By the~above, for neither of $a\in D$ and $i=1,\dots,n$,
the~inequalities $m(a)\ge \alpha_i$ and $\phi(a)(b)\ge\beta_i$ are
valid simultaneously. The~set
$$
B_i=\{a\in D\mid \phi(a)(b)\ge\beta_i\}=
\{a\in D\mid \phi(a)\ge \beta_i*\eta_{[L]}D'(b)\}
$$
is closed w.r.t.\ the~lower topology due to the~continuity of
$\phi$. It has an empty intersection with the~Scott closed set
$$
A_i=\{a\in D\mid m(a)\ge\alpha_i\}.
$$

By compactness, there is a~finite collection
$a_{i1},\dots,a_{ik_i}\in D$ such that the set
$$
\{a\in D\mid a_{ij}\le a\text{ for some }1\le j\le k_i\}
$$
contains $B_i$ and has an~empty intersection with $A_i$. Then
the~set
$$
V=
\{c\in \uM_{[L]}D\mid c\not\ge \alpha_i*\eta_{[L]}D(a_{ij})
\text{ for all }1\le i\le n,1\le j\le k_i\}
$$
is an~open neighborhood of $m$ in the~lower topology, and, if $c\in
V$, then $c(a)\not \ge \alpha_i$ whenever $\phi(a)(b)\ge\beta_i$,
$1\le i\le n$.

Therefore, if $c\in V$, then
$$
\sup\{c(a)*\phi(a)(b)\mid a\in D\}\le\gamma'\not \ge m'(b),
$$
hence
$
\usp(\phi)(c)(b)\not \ge m'(b),
$
and all preimages $\usp(\phi)^{-1}(\{m'\}\ups)$ are closed, which
implies the~required continuity of~$\usp(\phi)$.
\end{proof}

\begin{proposition}\label{stat.lin-sp}
Let $\phi$ be a~mapping $D\to\uM_{[L]}D'$. If (a) $\phi$ is Scott
continuous, or (b) $*$ is infinitely distributive w.r.t.\ infimum,
then the~mapping $\usp(\phi):\uM_{[L]}D\to\uM_{[L]}D'$ is linear.
\end{proposition}

\begin{proof}
Join preservation is due to Proposition~\ref{stat.scott-sp}.

Let a~mapping $\phi:D\to\uM_{[L]}D'$ be Scott continuous (a). Then:
\begin{gather*}
\usp(\phi)(\alpha\bd m)=
\usp(\phi)((\alpha * m)^u)\mathrel{\overset{\scriptstyle\rm Lemma~\ref{lem.post-u}}
{=\kern-2pt=\kern-2pt=\kern-2pt=\kern-2pt=\kern-2pt=}}
\usp(\phi)(\alpha * m)=
\\
\bigl(\psup_{a\in D}\alpha*m(a)*\phi(a)\bigr)^u
\mathrel{\overset{\scriptstyle\rm Lemma~\ref{lem.alpha-u}}
{=\kern-2pt=\kern-2pt=\kern-2pt=\kern-2pt=\kern-2pt=}}
\bigl(\alpha*(\psup_{a\in D}m(a)*\phi(a))^u\bigr)^u=
\alpha\bd\usp(\phi)(m).
\end{gather*}

Assume (b). Then:
\begin{gather*}
\usp(\phi)(\alpha\bd m)(b)= \usp(\phi)(\alpha* m)(b)=
\\
\inf\bigl\{\sup\{\alpha*m(a)*\phi(a)(b') \mid a\in D\}
\mid b'\in D',b'\ll b\bigr\}=
\\
\inf\bigl\{\alpha*\sup\{m(a)*\phi(a)(b') \mid a\in D\}
\mid b'\in D',b'\ll b\bigr\}=
\\
\alpha*\inf\bigl\{\sup\{m(a)*\phi(a)(b') \mid a\in D\}
\mid b'\in D',b'\ll b\bigr\}=
\\
\alpha\bd \usp(\phi)(m)(b),\text{ \ \ for all }m\in\uM_{[L]}D, b\in D'.
\end{gather*}
\end{proof}

\begin{remark}
In the~presence of (a) or (b), the~mapping $\usp(\phi)$ can be
characterized as the~least linear mapping $\Phi:\uM_{[L]}D\to
\uM_{[L]}D'$ such that $\Phi(\eta_{[L]}D(d))= \phi(d)$ for all
$d\in D$.
\end{remark}

\begin{remark}
\emph{All} statements in this section have straightforward
analogues for normalized predicates. The only significant
distinction is that, if a~mapping $\phi:D\to M_{[L]}D'$ satisfies
the~conditions that are analogous to ones of \ref{stat.lin-sp},
then the~mapping $sp(\phi):M_{[L]}D\to M_{[L]}D'$ is affine,
instead of linear. Proofs can be obtained
\emph{mutatis mutandis}, without any major changes.
\end{remark}

\section*{Epilogue}

We have shown that $L$-fuzzy strongest postcondition predicate
transformers are related to $L$-idempotent linear or affine
operators between continuous $L$-semimodules. Now it is possible to
study linear and affine approximations of predicate transformers
from above and from below. These approximations are related to
attempts to describe a~program behaviour in a~more economical way,
dropping less important details.

It has been observed, e.g., by
Doberkat~\cite{Dob:DemProdProbRel:01} that monads and Kleisli
composition arise in description of combining several programs into
a pipe and composing the~respective predicate transformers. While,
for probabilistic programs, these monads are based on
(sub)probability measures, for non-probabilistic fuzzy semantics we
propose to use monads of lattice-valued non-additive
measures~\cite{Nyk:CapLat:08}.

Treatment of $L$-fuzzy weakest precondition predicate transformers,
similar to a~proposed one for strongest precondition predicate
transformers, as well as a~demonstration that relations between
these classes can be properly expressed in terms of category
theory, will be the~topic of our future publications. In
particular, Galois connections \cite{Nyk:ContDualContLsem:12} will
be used to investigate compatibility of $L$-fuzzy knowledge and of
nondeterministic programs.


\begin{thebibliography}{99}

\bibitem{Akian:DensInvMeas:99} M. Akian, \emph{Densities of invariant measures and large
deviations}, Trans. Amer. Math. Soc. \textbf{351}:11 (1999)
4515--4543.

\bibitem{Berg:IntrManValFuzLogic:08}
M.~Bergmann, Introduction to Many-Valued and Fuzzy Logic:
Semantics, Algebras, and Derivation Systems, Cambridge University
Press, N.Y., 2008.

\bibitem{ChWu:ImplBasFuzzPredTr:05}
Y.~Chen, H.~Wu, \emph{Implication-based fuzzy predicate
transformers}, in: Fuzzy logic, soft computing and computational
intelligence: Eleventh International Fuzzy Systems Association,
Vol.I, Beijing, 2005, pp.~77--82.

\bibitem{ChWu:SemWlpSlp:11}
Y.~Chen, H.~Wu, \emph{The Semantics of wlp and slp of Fuzzy
Imperative Programming Languages}, in: Nonlinear Mathematics for
Uncertainty and its Applications, Vol. 100 in Advances in
Intelligent and Soft Computing, 2011, pp.~357--364.

\bibitem{DMR:LattValPredTrIntSys:10}
J.T.~Denniston, A.~Melton, S.E.~Rodabaugh, \emph{Lattice-valued
predicate transformers and interchange systems}, in: P. Cintula, E.
P.~Klement, L.N.~Stout, Abstracts of the 31th Linz Seminar
(February 2010), Universit\"atsdirecktion Johannes Kepler
Universit\"at (Linz, Austria), 31--40.

\bibitem{Dijk:Guarded:75}
E.W.~Dijkstra, \emph{Guarded commands, non-determinacy and formal
derivation of programs}, Comm. of the ACM {\bf 18}:8 (1975)
453--457.

\bibitem{Dob:DemProdProbRel:01}
E.E. Doberkat, \emph{Demonic product of probabilistic relations},
Technical Report 116, University of Dortmund.

\bibitem{DubPr:GradUncBip-MakSense:10} D.~Dubois,
H.~Prade, \emph{Gradualness, uncertainty and bipolarity: Making
sense of fuzzy sets}, Fuzzy Sets and Systems, publ. online 19 Nov
2010, \texttt{doi:10.1016/j.fss.2010.11.007}.

\bibitem{Eda:DomComput:97}~A. Edalat, \emph{Domains for computation in
mathematics, physics and exact real arithmetic}, Bull. Symb. Logic
{\bf 3}:4 (1997) 401--452.

\bibitem{Erne-ZdistrFunSp:98}
M.~Ern\'e, Z-distributive function spaces, preprint, 1998.

\bibitem{GHKLM:ContLattDom:03}
G. Gierz, K. H. Hofmann, K. Keimel, J. D. Lawson, M. Mislove, D. S.
Scott, Continuous Lattices and Domains, Cambridge University Press,
2003.

\bibitem{Haj:FuzNoncomLog:03}
P.~H\'ajek, \emph{Fuzzy logics with noncommutative conjuctions} J
Logic Computation {\bf 13}:4 (2003) 469--479.

\bibitem{HH:DualTheoQuantSem:98}
R.~Heckmann, M.~Huth,
\emph{A~duality theory for quantitative semantics}, in:
Proceedings of the 11th International Workshop on Computer Science
Logic, volume 1414 of Lecture Notes in Computer Science, Springer
Verlag, 1998, pp.~255--274.

\bibitem{ML:CWM:98}
S.~Mac Lane, Categories for the Working Mathematician. 2nd ed.
Springer, New York, 1998.

\bibitem{MinY:ReasProbProgr:03}
Mingsheng Ying, \emph{Reasoning about probabilistic sequential
programs in a probabilistic logic}, Acta Inf. {\bf 39}:5 (2003)
315--389.

\bibitem{MMS:ProbPredTrans:96}
C.~Morgan, A.~McIver, K.~Seidel, \emph{Probabilistic predicate
transformers}, ACM Trans. Program. Lang. Syst. {\bf 18}:3 (1996)
325--353.

\bibitem{Nyk:CapLat:08} O.~Nykyforchyn,
\emph{Capacities with values in compact Hausdorff lattices},
Appl. Cat. Struct. {\bf 15}:3 (2007) 243--257.

\bibitem{Nyk:AdMonCompLawsSem:11}
O.~Nykyforchyn, \emph{Adjoints and monads related to compact
lattices and compact Lawson idempotent semimodules}, Order {\bf
29}:1 (2012) 192--213.

\bibitem{Nyk:ContDualContLsem:12}
O.~Nykyforchyn, \emph{Continuous and dually continuous
$L$-semimodules}, Mat.Stud. {\bf 37}:1 (2012) 3--28.

\bibitem{Rosen:Quantales:90}
K.~Rosenthal, Quantales and Their Applications, Pitman Research
Notes in Mathematics Series 234, Longman Scientific \& Technical,
Wiley, Essex, England, New York, 1990.

\bibitem{Verw:RandUSC:97}
W.~Verwaat, \emph{Random upper semicontinuous functions and
extremal processes} in: Probability and Lattices. Eds: W.~Verwaat,
H.~Holverda. CWI tract, Centrum voor Wiskunde en Informatica, 1997,
pp.~1--56.

\end{thebibliography}
\end{document}